\input amstex\documentstyle{amsppt}  
\pagewidth{12.5cm}\pageheight{19cm}\magnification\magstep1
\topmatter
\title The quantum group $\dot U$ and flag manifolds over the
semifield $\bold Z$\endtitle
\author G. Lusztig\endauthor
\address{Department of Mathematics, M.I.T., Cambridge, MA 02139}\endaddress    
\thanks{Supported by NSF grant DMS-2153741}\endthanks
\endtopmatter   
\document

\define\ds{\dot s}

\define\frl{\forall}

\define\si{\sim}

\define\sqc{\sqcup}

\define\qua{\quad}

\define\lb{\linebreak}

\define\op{\oplus}
   
\define\part{\partial}
\define\emp{\emptyset}

\define\n{\notin}
\define\iy{\infty}
\define\m{\mapsto}
\define\do{\dots}

\define\sub{\subset}    

\define\T{\times}
\define\ti{\tilde}
\define\nl{\newline}
\redefine\i{^{-1}}
\define\fra{\frac}
\define\un{\underline}

\define\ot{\otimes}

\define\Hom{\text{\rm Hom}}

\redefine\b{\beta}

\define\g{\gamma}
\redefine\d{\delta}
\define\e{\epsilon}
\define\et{\eta}
\define\io{\iota}
\redefine\o{\omega}
\define\p{\pi}
\define\ph{\phi}

\redefine\t{\tau}
\define\th{\theta}
\define\k{\kappa}
\redefine\l{\lambda}
\define\z{\zeta}
\define\x{\xi}

\define\Om{\Omega}

\define\Th{\Theta}
\redefine\L{\Lambda}

\redefine\aa{\bold a}

\define\bof{\bold f}
\define\hh{\bold h}
\define\ii{\bold i}

\define\kk{\bold k}

\define\uu{\bold u}

\redefine\AA{\bold A}
\define\BB{\bold B}

\define\NN{\bold N}

\define\QQ{\bold Q}
\define\RR{\bold R}

\define\ZZ{\bold Z}

\define\cb{\Cal B}

\define\ce{\Cal E}
\define\cf{\Cal F}

\define\ci{\Cal I}

\define\ck{\Cal K}

\define\cu{\Cal U}

\define\ta{\ti a}

\define\tB{\ti B}

\define\tE{\ti E}
\define\tF{\ti F}

\define\tK{\ti K}

\define\sha{\sharp}

\define\dU{\dot U}
\define\dBB{\dot\BB}

\head Introduction\endhead
\subhead 0.1 \endsubhead
Let $\bof$ be the $+$ part of the Drinfeld-Jimbo quantized
enveloping algebra $U$ (over $\QQ(v)$)
attached to a root datum of simply laced
type and let $\dU$ be the modified form of $U$ (see
\cite{L93, \S23}). Let $\BB$ (resp. $\dBB$) be the canonical basis
of $\bof$ (resp. $\dU$) defined in \cite{L90a} (resp. in
\cite{L92}, see also \cite{L93, 25.2}). In
\cite{L90a} it was shown that $\BB$ is 
naturally parametrized by something which later \cite{L94} was
interpreted in terms of certain objects attached to the semifield
$\ZZ$. In this paper we want to find an analogous parametrization
for $\dBB$ (see 5.12), compatible with the involution $\o$ of
$\dU$ interchanging the $+$ part and the $-$ part (see 3.3)
which preserves $\dBB$ (see 3.8). In particular we show that
$\dBB$ is preserved by $\o$.

Let $X$ be the lattice of weights of our root datum and let
$X^+$ be the set of dominant weights in $X$. We will put $\dBB$
in bijection with
$\sqc_{\l\in X^+}(B_\l\T B_\l)$
where $B_\l$ is the canonical basis \cite{L90a} of
the simple finite dimensional $U$-module $\L_\l$ with
highest weight $\l$. We are reduced to finding a parametrization
of $B_\l$ in terms of objects attached to the semifield $\ZZ$.
Such a parametrization was given in \cite{L90a,\S8}, but this is
still not enough for our purpose.

We will give another parametrization of $B_\l$ based on the
following observation of \cite{L97}.
According to \cite{L90a,\S8},
the set $B_\l$ can be parametrized in two different ways:
by regarding $\L_\l$ as a highest weight module or as a lowest
weight module. The rather complicated combinatorics relating these
two parametrizations was shown in \cite{L97} to be expressible
in terms of a remarkable involution $\ph_\ZZ$ (defined in
\cite{L97}) of a certain object $\cb_\ZZ$ attached to the
semifield $\ZZ$.
We can parametrize $B_\l$ is in terms of this involution
$\ph_\ZZ$. More precisely, if $G$ is a reductive group
corresponding to the dual root datum, then $\cb_\ZZ$ is
a variant of the flag manifold of $G$
over the semifield $\ZZ$. Then $\cb_\ZZ$ has two remarkable
subsets $\cb_\NN^+,\cb_\NN^-$ (interchanged by
$\ph_\ZZ$) in which certain parameters in $\ZZ$ are assumed to be
in $\NN$. As in \cite{L97}, $\cb_\ZZ$ has
an action of the group $X$; this is the $\ZZ$-variant of the
conjugation action of a
maximal torus of $G$ on the flag manifold of $G$. Now, for
$\l\in X$ we can consider the intersection of $\cb_\NN^+$ with the
$\l$ translate (in the sense of the $X$-action)
of $\cb_\NN^-$. Although $\cb_\NN^+$ and $\cb_\NN^-$ are in
general infinite, we show that the intersection above is finite;
moreover, it is nonempty if and only if $\l$ is dominant,
 in which case it naturally parametrizes $B_\l$.
This parametrization is an immediate consequence of
\cite{L97, 4.9}, except for the nonemptiness criterion above.
(The proof of 4.9 in \cite{L97} contained some 
misprints and we reprove it in 5.5.)

The resulting parametrization of $\dBB$ has the property that
the action of $\o$ has a simple description in terms of the 
involution $\ph_\ZZ$.

Another parametrization of $B_\l$ (which again uses ideas in
\cite{L97}) is given in \cite{GS15}.

\subhead 0.2 \endsubhead
One can show that similar results hold when our root datum is not
assumed to be simply laced, by using folding to reduce to the
simply laced case (as was done for $\BB$ in \cite{L11} where again
the group $G$ corresponding to the dual root datum was used.)

\head 1. The flag manifold $\cb_K$ for a semifield $K$\endhead
\subhead 1.1 \endsubhead
We fix a finite set $I$ and a simply laced Cartan matrix $(a_{ij})$
indexed by $I\T I$, that is a symmetric, positive definite matrix
with integer entries such that $a_{ii}=2$ for $i\in I$,
$a_{ij}\in\{0,-1\}$ for $i\ne j$ in $I$. We denote the obvious
$\ZZ$-basis of $\ZZ[I]$ as $\{i';i\in I\}$. For $i\in I$ we define
a homomorphism $s_i:\ZZ[I]@>>>\ZZ[I]$ by
$j'\m j'-a_{ij}i'$ ($j\in I$). Let $W$ be the subgroup of the
automorphism group of $\ZZ[I]$ generated by $\{s_i;i\in I\}$. This
is a finite Coxeter group with length function $w\m|w|$.

Let $w_0$ be the unique element of $W$ with $|w_0|$ maximal and let
$\nu=|w_0|$. 

For $i\in I$ we define $i^!\in I$ by $s_{i^!}=w_0s_iw_0$; then
$i\m i^!$ is an involution of $I$. 

Let $\ci$ be the set of sequences $\ii=(i_1,i_2,\do,i_\nu)$ such
that $s_{i_1}s_{i_2}\do s_{i_\nu}=w_0$. For example, if
$I=\{i,j\},a_{ij}=-1$, we have $\ci=\{(i,j,i),(j,i,j)\}$.

\subhead 1.2 \endsubhead
Let $K$ be a semifield. Let $(\ii,\aa),(\ii',\aa')$ in
$\ci\T K^\nu$ be given by
$$\ii=(i_1,i_2,\do,i_\nu), \ii'=(i'_1,i'_2,\do,i'_\nu),$$
$$\aa=(a_1,a_2,\do,a_\nu), \aa'=(a'_1,a'_2,\do,a'_\nu).$$
We say that $(\ii,\aa),(\ii',\aa')$ are adjacent if one of (i),(ii)
below holds.

(i) For some $l\in[1,\nu-3]$ we have $i'_k=i_k$ for
$k\n\{l,l+1,l+2\}$ and $(i_{l},i_{l+1},i_{l+2})=(i,j,i)$,
$(i'_{l},i'_{l+1},i'_{l+2})=(j,i,j)$ where $i,j$ in $I$ satisfy
$a_{ij}=-1$; moreover, $(a_l,a_{l+1},a_{l+2})=(a,b,c)$,
$(a'_l,a'_{l+1},a'_{l+2})=(a',b',c')$ where
$$a'=bc/(a+c),b'=a+c,c'=ab/(a+c)$$,
or equivalently
$$a=b'c'/(a'+c'),b=a'+c',c=a'b'/(a'+c');$$

(ii) for some $l\in[1,\nu-2]$ we have $i'_k=i_k$ for
$k\n\{l,l+1\}$ and $(i_{l},i_{l+1})=(i,j)$,
$(i'_{l},i'_{l+1})=(j,i)$ where $i,j$ in $I$ satisfy $a_{ij}=0$;
moreover, $a'_l=a_{l+1}$, $a'_{l+1}=a_l$.

Let $\cu_K$ be the set of equivalence classes on $\ci\T K^\nu$ for
the equivalence relation on $\ci\T K^\nu$ generated by the
adjacency relation. 

We shall sometime denote $(\ii,\aa)\in\ci\T K^\nu$, (or its
equivalence class) as \lb
$i_1^{a_1}i_2^{a_2}\do i_\nu^{a_\nu}$, where
$\ii=(i_1,i_2,\do,i_\nu),\aa=(a_1,a_2,\do,a_\nu)$.

The assignment
$$i_1^{a_1}i_2^{a_2}\do i_\nu^{a_\nu}@>>>
i_\nu^{a_\nu}i_{\nu-1}^{a_{\nu-1}}\do i_1^{a_1}$$
defines an involution $A\m A^*$ of $\cu_K$.

\subhead 1.2\endsubhead
For $\ii=(i_1,i_2,\do,i_\nu)\in\ci$ and $k\in[1,\nu]$, we have
$$s_{i_1}s_{i_2}\do s_{i_{k-1}}(i'_k)=\sum_{h\in I}r_{h,k}h'$$
(in $\ZZ[I]$) where $r_{h,k}\in\NN$.
For $\aa=(a_1,a_2,\do,a_\nu)\in K^\nu, h\in I$ we set
$$||\ii,\aa||_h=\prod_{k\in[1,\nu]}a_k^{r_{h,k}}\in K.$$
We show:

(a) If $(\ii,\aa),(\ii',\aa')$ in $\ci\T K^\nu$ are adjacent, then
$||\ii,\aa||_h=||\ii',\aa'||_h$.
\nl
Let $r'_{h,k}\in\NN$ be defined in terms of $\ii'$ in the same way
as $r_{h,k}$ was defined in terms of $\ii$. In case 1.1(i) we must
show:
$$a_l^{r_{h,l}}a_{l+1}^{r_{h,l+1}}a_{l+2}^{r_{h,l+2}}
=a'_l{}^{r'_{h,l}}a'_{l+1}{}^{r'_{h,l+1}}a'_{l+2}{}^{r'_{h,l+2}}.
\tag b$$

For some $w\in W$ we have

$w(i')=\sum_{h_1\in I}r_{h_1,l}h'_1$,

$ws_i(j')=w(i')+w(j')=\sum_{h_1\in I}r_{h_1,l+1}h'_1$,

$ws_is_j(i')=w(j')=\sum_{h_1\in I}r_{h_1,l+2}h'_1$,

$w(j')=\sum_{h_1\in I}r'_{h_1,l}h'_1$,

$ws_j(i')=w(i')+w(j')=\sum_{h_1\in I}r'_{h_1,l+1}h'_1$,

$ws_js_i(j')=w(i')=\sum_{h_1\in I}r'_{h_1,l+2}h'_1$.

Thus,

$r_{h,l}=r'_{h,l+2}$, $r_{h,l+1}=r'_{h,l+1}=r_{h,l}+r_{h,l+2}$,
$r_{h,l+2}=r'_{h,l}$.

Thus, (b) is equivalent to

$$(a_la_{l+1})^{r_{h,l}}(a_{l+1}a_{l+2})^{r_{h,l+2}}
=(a'_la'_{l+1})^{r_{h,l}}(a'_{l+1}a'_{l+2})^{r_{h,l+2}}.$$

This follows from 
$$a_la_{l+1}=a'_{l+1}a'_{l+2},\qua a'_la'_{l+1}=a_{l+1}a_{l+2}.$$
This proves (a) in case 1.1(i).

In case 1.1(ii) we must show
$$a_l^{r_{h,l}}a_{l+1}^{r_{h,l+1}}=
a'_l{}^{r'_{h,l}}a'_{l+1}{}^{r'_{h,l+1}}.\tag c$$
For some $w\in W$ we have

$w(i')=\sum_{h_1\in I}r_{h_1,l}h'_1$,

$ws_i(j')=w(j')=\sum_{h_1\in I}r_{h_1,l+1}h'_1$,

$w(j')=\sum_{h_1\in I}r'_{h_1,l}h'_1$,

$ws_j(i')=w(i')=\sum_{h_1\in I}r'_{h_1,l+1}h'_1$.

Thus, $r_{h,l}=r'_{h,l+1}$, $r_{h,l+1}=r'_{h,l}$ so that

$r_{h,l}=r'_{h,l+1}$, $r_{h,l+1}=r'_{h,l}$.

Thus, (c) is equivalent to
$$a_l^{r_{h,l}}a_{l+1}^{r_{h,l+1}}=
a'_l{}^{r_{h,l+1}}a'_{l+1}{}^{r_{h,l}}$$

and this follows from $a'_l=a_{l+1}$, $a'_{l+1}=a_l$. This proves
(a) in case 1.1(ii).

In view of (a) we can define $||A||_h\in K$ for any
$A\in\cu_K,h\in I$ to be $||\ii,\aa||_h$ for any $(\ii,\aa)$ in $A$.

\subhead 1.3 \endsubhead
Let $\ii\in\ci$. We define a map $f_\ii:K^\nu@>>>\cu_K$ by sending
$\aa\in K^\nu$ to the equivalence class of $(\ii,\aa)$. We show:

(a) $f_\ii$ is a bijection.
\nl
Let $A\in\cu_K$. If $(\ii',\aa')\in A$ then using Matsumoto's
theorem \cite{L03, 1.9} for $w_0$, we can find a sequence
$(\ii^1,\aa^1),\do,(\ii^n,\aa^n)$ in $\ci\T K^\nu$ in which any
two consecutive terms are adjacent such that
$(\ii^n,\aa^n)=(\ii',\aa')$, $\ii^1=\ii$.
We have $(\ii,\aa^1)\in A$. Thus $f_\ii$ is surjective.

Next we prove that $f_\ii$ is injective.
Assume first that $K$ is contained on the multiplicative group of a
field $\kk$ of characteristic zero with $+,\T$ induced from that
of $\kk$. Then $f_\ii$ is injective by 4.2(a).

Next we consider a general $K$.
Assume that $\aa=(a_1,\do,a_\nu)\in K^\nu,
\aa'=(a'_1,\do,a'_\nu$ in $K^\nu$ are such that
$f_\ii(\aa)=f_\ii(\aa')$. Then we can find
$$(\ii^1,\aa^1),\do,(\ii^n,\aa^n)$$ in $\ci\T K^\nu$ in which any two
consecutive terms are adjacent and $(\ii^1,\aa^1)=(\ii,\aa)$,
$(\ii^n,\aa^n)=(\ii,\aa')$.
By \cite{BFZ, 2.1.6} we can find a semifield $\tK$ as in the first
part of the proof and a homomorphism od semifields $z:\tK@>>>K$
such that $$a_1=z(\ta_1),\do,a_\nu=z(\ta_\nu)$$ for some
$$\ti\aa=(\ta_1,\do,\ta_\nu)\in\tK^\nu.$$
We define $$\ti\aa^1,\ti\aa^2,\do,\ti\aa^\nu$$ by the condition that
any two consecutive terms of
$$(\ii^1,\ti\aa^1),(\ii^2,\ti\aa^2),\do,(\ii^n,\ti\aa^n)$$ are
adjacent (in $\ci\T\tK^\nu$). Now $z$ takes $\ti\aa^1$ to $\aa^1$
and then it automatically takes $\ti\aa^2,\do,\ti\aa^n$ to
$\aa^2,\do,\aa^n$. We have $\ii^1=\ii^n=\ii$ and the injectivity
of our map (for $\tK$) implies that $\ti\aa^1=\ti\aa^n$. Applying
$z$ we see that $\aa^1=\aa^n$ that is $\aa=\aa'$. This proves
injectivity of $f_\ii$. This proves (a).

\subhead 1.4 \endsubhead
Let $i\in I,c\in K$.
If $A\in\cu_K$ we can find $(\ii,\aa)\in A$ such that the first
term of $\ii$ is $i$. (We use Matsumoto's theorem \cite{L03, 1.9}
for $w_0$.)
We set ${}_c\aa=(a_1c,a_2,\do,a_\nu)$ where $\aa=(a_1,a_2,\do,a_\nu)$.
Assume that  $(\ii',\aa')\in A$ is also such that the first term of
$\ii'$ is $i$. Then ${}_c\aa'\in K^\nu$ is defined. We show that

(a) $(\ii,{}_c\aa),(\ii',{}_c\aa')$ are equivalent.
\nl
Using Matsumoto's theorem for $s_iw_0$, we see that we can find a
sequence
$$(\ii^1,\aa^1),\do,(\ii^n,\aa^n)$$ in $\ci\T K^\nu$ in which any two
consecutive terms are adjacent as in 1.1(i),(ii) with $l\ge2$ such
that $$(\ii^1,\aa^1)=(\ii,\aa), \ii^n=\ii'$$ (in particular each
$\ii^1,\ii^2,\do$ starts with $i$).
Then $(\ii',\aa^n),(\ii',\aa')$ are both in $A$ hence by 1.3(a)
we must have $\aa^n=\aa'$. Now any two consecutive terms of 
$$(\ii^1,{}_c\aa^1),\do,(\ii^n,{}_c\aa^n)$$ are adjacent and we have
$$(\ii^1,{}_c\aa^1)=(\ii,{}_c\aa),(\ii^n,{}_c\aa^n)=
(\ii',{}_c\aa').$$
This proves (a).

We see that $(\ii,\aa)\m(\ii,\aa_c)$ defines a map (in fact a
bijection) $T_{i,c}:\cu_K@>>>\cu_K$.

For $c,c'$ in $K$ we have $T_{i,c}T_{i,c'}=T_{i,cc'}$.

\subhead 1.5 \endsubhead
Let $i\in I,c\in K,A\in\cu_K$

We can find $(\ii,\aa)\in A$ such that the first term of $\ii$ is
$i$. Define $r_{h,k}$ in terms of $\ii$ as in 1.2. We have
$r_{h,1}=\d_{h,i}1$ (this is $1$ if $h=i$ and is $0$ if $h\ne i$.
We have
$$||\ii,\aa||_h=a_1^{\d_{i,h}}\prod_{k\in[2,\nu]}a_k^{r_{h,k}}\in K.$$
$$||\ii,{}_c\aa||_h=(ca_1)^{\d_{i,h}}\prod_{k\in[2,\nu]}a_k^{r_{h,k}}
=c^{\d_{i,h}}||\ii,\aa||_h$$
Thus we have

$||T_{i,c}A||_h=c^{\d_{i,h}}||A||_h$.

\subhead 1.6\endsubhead
We regard $K^I$ as a group (the product of $(p_i)$ and $(p'_i)$
is $(p_ip'_i)$).

For $p=(p_i)_{i\in I}\in K^I$ there is well defined bijection
$S_p:\cu_K@>>>\cu_K$ given by
$$i_1^{a_1}i_2^{a_2}\do i_\nu^{a_\nu}\m
i_1^{a_1p_{i_1}}i_2^{a_2p_{i_2}}\do i_\nu^{a_\nu p_{i_\nu}}$$
for $(\ii,\aa)\in\ci\T K^\nu$, see \cite{L21, no.8}.

This defines an action of the group $K^I$ on $\cu_K$.

For $p\in K^I,i\in I,c\in K$ we have

(a) $T_{i,c}S_p=S_pT_{i,c}$.
\nl
There is a well defined involution $\io:\cu_K@>>>\cu_K$ such that
$$i_1^{a_1}i_2^{a_2}\do i_\nu^{a_\nu}\m
(i_1^!)^{a_1}(i_2^!)^{a_2}\do(i_\nu^!)^{a_\nu}$$
for $(\ii,\aa)\in\ci\T K^\nu$. 
For $p\in K^I$ we have $\io S_p=S_{p^!}\io$
where $p^!\in K^I$ is given by $(p^!)_i=p_{i^!}$.

For $i\in I,c\in K$ we have

(b) $T_{i,c}\io=\io T_{i^!,c}$

We show:

(c) The action of $K^I$ on $\cu_K$ described above is free.
\nl
We must show that if $p=(p_i)\in K^I$ and $A\in\cu_K$ are such that
$S_pA=A$ then $p_i=1$ for all $i$. Let $(\ii,\aa)\in A$ with
$\ii=(i_1,i_2,\do,i_\nu),\aa=(a_1,a_2,\do,a_\nu)$. Let
$\aa'=(a_1p_{i_1},a_2p_{i_2},\do,a_\nu p_{i_\nu})\in K^\nu$.
By assumption we have $(\ii,\aa')\in A$. Using 1.3(a) we deduce
that $\aa'=\aa$. Thus $a_k=a_kp_{i_k}$ for $k=1,\do,\nu$ so that
$p_{i_k}=1$ for $k=1,\do,\nu$. For any $i\in I$ we can find $k$
such that $i_k=i$. It follows that $p_i=1$ for $i\in I$. This
proves (c).

\subhead 1.7 \endsubhead
For a semifield $K$, let $\ph_K:\cu_K@>>>\cu_K$ be the bijection
defined in 4.6, see also \cite{L21, no.11}. We have $\ph_K^2=1$.

For example, if $I=\{i\}$ we have $\ph_K(i^a)=i^{1/a}$; if
$I=\{i,j\}$ with $a_{ij}=1$, we have
$$\ph_K(i^aj^bi^c)
=i^{a/c(a+c)}j^{(a+c)/ab}i^{1/(a+c)}=j^{c/ab}i^{1/c}j^{1/b}.$$
We have

(a) $\io\ph_K=\ph_K\io$.
\nl
For $p\in K^I$ we have

(b) $S_p\ph_K=\ph_KS_{p\i}$;
\nl
hence  $S_p\ph_K:\cu_K@>>>\cu_K$ has square $1$. 

For $i\in I,c\in K$ we have
$$T_{i,c}\ph_K=\ph_KT_{i^!,c\i}.\tag c$$
Now (b),(c) can be viewed as sequences of equalities between certain
rational functions which are quotients of nonzero polynomials in
several variables with all coefficients in $\NN$ (after
substituting
elements of $K$ for the variables). It is enough to prove these
equalities in the case where $K=\RR_{>0}$; in that case (b) is
proved in \cite{L19, 4.3(d)} and (c) is proved in \cite{L21, 10(a)}.

\subhead 1.8\endsubhead
Let $\cb_{K}=\{(A,A')\in\cu_K\T\cu_K;\ph_K(A)=A'\}$ be the graph of  
$\ph_K$. We define an involution $\un\ph_K:\cb_{K}@>>>\cb_{K}$ by
$\un\ph_K(A,A')=(A',A)$. (We use that $\ph_K^2=1$.)
We define an involution $\un\io:\cb_{K}@>>>\cb_{K}$ by
$\un\io(A,A')=(\io(A),\io(A'))$. (We use 1.7(a).)
Now the group $K^I$ acts on $\cb_{K}$ by $p:B\m\un S_p(B)$ where
$\un S_p(A,A')=(S_pA,S_{p\i}A')$. (We use 1.7(b)).

For $p\in K^I$ we have
$\un S_p\un\ph_K=\un\ph_K\un S_{p\i}$.

\subhead 1.9\endsubhead
Let $A\in\cu_K$ and let $(\ii,\aa)\in A$ with
$\ii=(i_1,i_2,\do,i_\nu)$, $\aa=(a_1,a_2,\do,a_\nu)$.
Let $\ii'=(i_\nu,\do,i_2,i_1)\in\ci$. Define
$\aa'=(a'_1,a'_2,\do,a'_\nu)$ by $(\ii',\aa')\in\ph_K(A)$.

\proclaim{Lemma 1.10}We have
$a'_\nu=(\sum_{k\in[1,\nu];i_k=i_1}a_k)\i$.
\endproclaim
Using \cite{BFZ, 2.1.6} we see that we can reduce the general case
to the case where $K$ is as in 4.2. In that case the result follows
from 4.5. (More precisely 4.5 gives the analogous result for
$\ph'_K$ in 4.6 instead of 4.5, but the case of $\ph_K$ is then
a consequence.)

\head 2. The subsets $\cb^+_\NN,\cb^-_\NN$ of $\cb_\ZZ$\endhead
\subhead 2.1\endsubhead
In this section we assume that $K=\ZZ$ with the usual semifield
structure: the sum of $a,b$ is taken to be $\min(a,b)$; the product
of $a,b$ is taken to be $a+b$. Then $\cu_\ZZ$ is defined as in 1.1.
Note that the subset $\ci\T\NN^\nu$ of $\ci\T\ZZ^\nu$ is a union of
equivalence classes for $\si$. (We use that, if $a,c\in\NN$, then
$a-\min(a,c)\in\NN$.) The set of these equivalence classes is a
subset $\cu_\NN$ of $\cu_\ZZ$.

For $i\in I$ there is a well defined map $g'_i:\cu_\NN@>>>\NN$ such
that $i_1^{a_1}i_2^{a_2}\do i_\nu^{a_\nu}\m a_\nu$ whenever
$(\ii,\aa)\in\ci\T\NN^\nu$ satisfies $i_\nu=i$.

For any $p=(p_i)\in\NN^I$ we define
$$\cu_{\NN,p}=\{x\in\cu_\NN;g'_i(x)\le p_i\qua\frl i\in I\}.$$
Note that $\io:\cu_\ZZ@>>>\cu_\ZZ$ restricts to a bijection
$$\cu_{\NN,p}@>\si>>\cu_{\NN,p^!}$$
where $p^!\in\NN^I$ is as in 1.6.

\subhead 2.2\endsubhead
For $p=(p_i)\in\NN^I$ we set
$\hh_p=i_1^{p_{i_1}}i_2^{p_{i_2}}\do i_\nu^{p_{i_\nu}}\in\cu_\NN$
where $\ii\in\ci$. From the definition we see that this is
independent of the choice of $\ii$ hence is well defined.

Note that $\hh_0\in\cu_{\NN,p},\hh_p\in\cu_{\NN,p}$. We have
$$\hh_p=S(p)(\hh_0),\tag a$$
$$\io(\hh_p)=\hh_{p^!},\tag b$$
From \cite{L97, 2.9(a), 3.9} we have
$$\ph_\ZZ(\hh_0)=\hh_0.\tag c$$
In fact, by \cite{L97, 2.9}, $\ph_\ZZ$ is the unique bijection
$\cu_\ZZ@>>>\cu_\ZZ$ satisfying (c) and 1.7(c) (with $K=\ZZ$).

\subhead 2.3\endsubhead
Let $p=(p_i)\in\NN^I$. By definition, for $h\in I$ we have
$$||\hh_p||_h=\sum_{k\in[1,\nu]}r_{h,k}p_{i_k}\in\NN$$
where $\ii=(i_1,i_2,\do,i_\nu)\in\ci$ and for $k\in[1,\nu]$,
$r_{h,k}\in\NN$ are defined by
$$s_{i_1}s_{i_2}\do s_{i_{k-1}}(i'_k)=\sum_{h\in I}r_{h,k}h'$$
(in $\ZZ[I]$). We state the following result.
$$||\hh_p||_h=\sum_{i\in I}(p_i+(p^!)_i)b_{ih}\in\NN\tag a$$
where $(b_{ij})$ is the inverse of the Cartan matrix $(a_{ij})$.
It is known that $b_{ij}\in\QQ_{\ge0}$ (see for example \cite{LT}).
Now (a) can be verified case by case.

Assume for example that
$I=\{i\}$. We have $||h_p||_i=p_i$,
$b_{ii}=1/2$, $(p_i+(p^!)_i)b_{ii}=p_i$ hence (a) holds.

Assume now that $I=\{i,j\}$ with $a_{ij}=-1$.
Let $\ii=(i,j,i)$. We have $\ii\in\ci$.
The corresponding sequence $\sum_{h\in I}r_{h,k}h'$, ($k=1,2,3$)
is $i',i'+j',j'$ hence
$$\sum_{h\in I}||\hh_p||_hh'
=p_ii'+p_j(i'+j')+p_ij'=(p_i+p_j)(i'+j').$$

We have $b_{ii}=b_{jj}=2/3,b_{ij}=b_{ji}=1/3$. Hence 
$$\align&\sum_{h\in I}(p_i+(p^!)_i)b_{ih}h'+
\sum_{h\in I}(p_j+(p^!)_j)b_{jh}h'\\&=
(p_i+p_j)((2/3)i'+(1/3)j')+(p_j+p_i)((1/3)i'+(2/3)j')
=(p_i+p_j)(i'+j').\endalign$$

Thus, (a) holds.

Assume now that $I=\{0,c,d,e\}$ with

$a_{0c}=a_{c0}=a_{0d}=a_{d0}=a_{0e}=a_{e0}=-1$,

$a_{cd}=a_{dc}=a_{ce}=a_{ec}=a_{de}=a_{ed}=0$.

Let
$\ii=(c,d,e,0,c,d,e,0,c,d,e,0)$. We have $\ii\in\ci$.

The corresponding sequence $\sum_{h\in I}r_{h,k}h'$ ($k=1,2,\do,12$)
is
$$c',d',e',0'c'd'e',0'd'e',0'c'e',0'c'd',0'0'c'd'e',0'c',0'd',0'e',
0'$$

where we omit $+$ signs (for example we write $0'c'$ for $0'+c'$.)

Thus we have
$$\align&\sum_{h\in I}||\hh_p||_hh'=
p_cc'+p_dd'+p_ee'+p_0(0'+c'+d+e')+p_c(0'+d'+e')\\&
+p_d(0'+c'+e')+p_e(0'+c'+d')+p_0(0'+0'+c'+d'+e')\\&
+p_c(0'+c')+p_d(0'+d')+p_e(0'+e')+p_00'+\\&
(2p_c+p_d+p_e+2p_0)c'+(p_c+2p_d+p_e+2p_0)d'\\&
+(p_c+p_d+2p_e+2p_0)e'+(2p_c+2p_d+2p_e+4p_0)0'.\endalign$$

We have

$b_{00}=2,b_{cc}=b_{dd}=b_{ee}=1$,

$b_{0c}=b_{c0}=b_{0d}=b_{d0}=b_{0e}=b_{e0}=1$,

$b_{cd}=b_{dc}=b_{ce}=b_{ec}=b_{de}=b_{ed}=1/2$.

Moreover, $p=p^!$. We see that (a) holds.

\subhead 2.4\endsubhead
We define two subsets of $\cb_{\ZZ}$ (see 1.8) by
$\cb^+_\NN=\{(A,A')\in\cb_\ZZ;A\in\cu_\NN\}$,
$\cb^-_\NN=\{(A,A')\in\cb_\ZZ;A'\in\cu_\NN\}$.

For any $p\in\ZZ^I$ we define
$$\cb_{\NN,p}=\cb^+_\NN\cap\un S_p(\cb^-_\NN)\sub\cb_{\ZZ}.$$

Clearly, $(A,A')\m A$ is a bijection
$$\cb_{\NN,p}@>\si>>\{A\in\cu_\NN;S_p\ph_\ZZ(A)\in\cu_\NN\}.\tag a$$
Note that if $p\in\NN^I$ then $(\hh_0,\hh_p)\in\cb_{\NN,p}$. (We use
2.2(a),(c).) In particular, in this case we have
$\cb_{\NN,p}\ne\emp$.

\subhead 2.5\endsubhead
We show that, conversely,

(a) if $p\in\ZZ^I$ and $\cb_{\NN,p}\ne\emp$, then $p\in\NN^I$.
\nl
Assume for example that $I=\{i,j\}$ with $a_{ij}=-1$.
and that $A\in\cu_\NN$ satisfies $S_p\ph_\ZZ(A)\in\cu_\NN$.
We show that $p\in\NN^I$.
We have $i^aj^bi^c\in A$ with $a,b,c$ in $\NN$ and
$j^{c-a-b+p_j}i^{-c+p_i}j^{-b+p_j}\in S_p\ph_\ZZ(A)$ with
$c-a-b+p_j\in\NN,-c+p_i\in\NN,-b+p_j\in\NN$. Then $p_i\ge c\ge0$,
$p_j\ge b\ge0$, so that indeed $p\in\NN^I$.

We now consider the general case.
Assume that $A\in\cu_\NN$ satisfies $S_p\ph_\ZZ(A)\in\cu_\NN$.
Let $i\in I$. We can find $(\ii,\aa)\in A$ with
$\ii=(i_1,i_2,\do,i_\nu)$, $\aa=(a_1,a_2,\do,a_\nu)\in\NN^\nu$
such that $i_\nu=i$.
Let $\ii'=(i_\nu,\do,i_2,i_1)\in\ci$. Define
$\aa'=(a'_1,a'_2,\do,a'_\nu)$ by $(\ii',\aa')\in\ph_K(A)$.
By 1.10 we have
$a'_\nu=-\min_{k\in[1,\nu];i_k=i_1}a_k$. In particular we have
$a'_\nu\le0$.
By assumption we have $p_i+a'_\nu\ge0$ hence $p_i\ge-a'_\nu\ge0$.
This proves (a).

\head 3. Canonical bases\endhead
\subhead 3.1 \endsubhead
We fix a simply laced root datum; this consists of a finite set $I$,
two finitely generated free abelian groups $Y,X$; a perfect pairing
$(,):Y\T X@>>>\ZZ$; an imbedding $I\sub X,i\m i'$ and an imbedding
$I\sub Y,i\m i$.
It is assumed that $I$ is as in 1.1 and $(i,j')=a_{ij}$ where
$a_{ij}$ is as in 1.1. We identify $\ZZ[I]$ in 1.1 with a subgroup
of $X$ by
$i'\m i'$. The action of $s_i$ on $\ZZ[I]$ extends to an action on
$X$ by $s_i:x\m x-(i,x)i'$. Here $i\in I$. Thus the $W$-action of
on $\ZZ[I]$ extends to a $W$-action on $X$.
For
$\l\in X$ we set $\l^!=-w_0(\l)$; then $\l\m\l^!$ is an involution
of $X$.

Let $X^+=\{\l\in X;(i,\l)\in\NN\qua\frl i\in I\}$. Note that $X^+$
is stable under $\l\m\l^!$.

For $\l,\l'$ in $X$ we write $\l'\ge\l$ if $\l'-\l\in\sum_{i\in I}\NN i'$
and $\l'>\l$ if $\l'\ge\l$ and $\l'\ne\l$.

\subhead 3.2\endsubhead
Let $v$ be an indeterminate. Let $\bof$ be the associative algebra with
$1$ over $\QQ(v)$ with generators $\{\th_i;i\in I\}$ associated to the
matrix $(a_{ij})$ in \cite{L94,1.2.5}. This can be identified with the
$+$ part of the algebra $U$ (see below) attached to the root datum. 

There is a unique algebra antiautomorphism $\bof@>>>\bof$ ($x\m x^*$)
such that $\th_i^*=\th_i$ for all $i\in I$. It has square $1$.

Let $U$ be the Drinfeld-Jimbo quantized enveloping algebra attached
to the root datum.
This is an associative algebra with $1$ over $\QQ(v)$.
As a vector space, $U$ can be identified with
$\op_{\g\in Y}\bof\ot\bof$ in two different ways: one by
$(x\ot x')_y\m x^+\ck_y x'{}^-$ and one by
$(x\ot x')_y\m x^-\ck_y x'{}^+$; here $\ck_0$ is the unit element.
The map $\bof@>>>U,x\m x^-$ and the map $\bof@>>>U,x\m x^+$ are
imbeddings of algebras with $1$.

\subhead 3.3\endsubhead
Let $\dU$ be the modified form of $U$, see \cite{L93, \S23}. This
is an associative algebra (without $1$ in general) over $\QQ(v)$.
In type $A$ it was defined in \cite{BLM}; the definition in the
general case is the same. As a vector space, $\dU$ can be
identified with
$\op_{\l\in X}\bof\ot\bof$ in two different ways: one by
$(x\ot x')_\l\m x^+1_\l x'{}^-$ and one by
$(x\ot x')_\l\m x^-1_\l x'{}^+$.

There is a unique vector space isomorphism $\sha:\dU@>>>\dU$ such that
$$\sha(x^+1_\l x'{}^-)=(x'{}^*)^+1_\l(x^*)^-$$ for $x,x'$ in $\bof$,
$\l\in X$;
we have also
$$\sha(x^-1_\l x'{}^+)=(x'{}^*)^-1_\l(x^*)^-$$ for $x,x'$ in $\bof$, $\l\in X$; hence $\sha^2=1$. Moreover, $\sha$ is an algebra antiautomorphism.

There is a unique vector space isomorphism $\o:\dU@>>>\dU$ such that
$$\o(x^+1_\l x'{}^-)=x^-1_{-\l}x'{}^+$$ for $x,x'$ in $\bof$, $\l\in X$; we have
also $\o(x^-1_\l x'{}^+)=x^+1_{-\l}x'{}^-$ for $x,x'$ in $\bof$, hence $\o^2=1$.
Moreover $\o$ is an algebra automorphism satisfying $\o\sha=\sha\o$.

\subhead 3.4\endsubhead
For $\l\in X^+$ let $\L_\l$ be the simple $U$-module defined in
\cite{L93, 3.5.6}. We shall regard $\L_\l$ also as (unital)
$\dU$-module as in \cite{L93, 23.1.4}. We have $\dim\L_\l<\iy$. Let
$\et_\l\in\L_\l$ be as in \cite{L93, 3.5.7}. We have
$1_\l\et_\l=\et_\l$.

Let $(,)_\l$ be the symmetric bilinear form $\L_\l\T\L_\l@>>>\QQ(v)$
defined in \cite{L93, 19.1.2}. Recall that $(\et_\l,\et_\l)_\l=1$.

Let $\dU[\ge\l]$ (resp. $\dU[>\l]$) be the set of all $u\in\dU$ such that
the following condition holds.

For any $\l'\in X^+$ such that $u$ acts on $\L_{\l'}$ by a nonzero map we
have $\l'\ge\l$ (resp.$\l'>\l$).
\nl
Clearly, $\dU[\ge\l]$ and $\dU[>\l]$ are two-sided ideals of $\dU$.

\subhead 3.5\endsubhead
Let $\BB$ be the canonical basis of $\bof$ (see \cite{L90a}, \cite{L93}).
By \cite{L90a, 3.3},

for $b\in\BB$ we have $b^*\in\BB$ and $b\m b^*$ is a bijection
$\BB@>\si>>\BB$.

If $b\in\BB$ then there is a well defined element
$wt(b)\in\sum_{h\in I}\NN h'\in X$
such that the following holds: $b$ is $\QQ(v)$-linear combination of
elements $\th_{i_1}\th_{i_2}\do\th_{i_n}$ where $i_1,i_2,\do,i_n$
in $I$ satisfy $i'_1+i'_2+\do+i'_n=wt(b)$.

\subhead 3.6\endsubhead
Let $\l\in X^+$. By \cite{L90a},

there is a unique $\QQ(v)$-basis $B_\l$ of $\L_\l$ and a unique
subset $\BB(\l)$ of $\BB$ such that $b\m b^-\et_\l$ maps
$\BB-\BB(\l)$ to $0$ and restricts to a bijection
$\BB(\l)@>\si>>B_\l$.

Let $\x_\l$ be the unique element in $B_\l$ such that
$1_{-\l^!}\x_\l=\x_\l$.

By \cite{L93, \S21},

(a) there is a unique vector space isomorphism
$\t:\L_\l@>>>\L_{\l^!}$ such that $\t(ux)=\o(u)\t(x)$ for
$u\in\dU,x\in\L_\l$ and $\t(\et_\l)=\x_{\l^!}$. It satisfies
$\t(B_\l)=B_{\l^!}$ and $\t(\x_\l)=\et_{\l^!}$.
\nl
We see that there is a unique bijection
$\k:\BB(\l)@>\si>>\BB(\l^!)$ such that  
$\t(b^-\et_\l)=\k(b)^-\et_{\l^!}$ for $b\in\BB(\l)$.
For $b\in\BB(\l)$ we have

(b) $b^+\x_{\l^!}=\k(b)^-\et_{\l^!}$.
\nl
This follows from $\t(b^-\et_\l)=b^+\x_{\l^!}$.

\subhead 3.7 \endsubhead 
From \cite{L93, 19.1.4} for any $b\in\BB(\l)$ we have
$$(b^-\et_\l,b^-\et_\l)_\l\in 1+v\i\QQ[[v\i]].$$   
   
\subhead 3.8 \endsubhead
Let $\dBB$ be the canonical basis of $\dU$ defined in
\cite{L92}, see also  \cite{L93, 25.2}. 
By \cite{L93, 26.3.2},

(a) for $\b\in\dBB$
we have $\o\sha(\b)\in\pm\dBB$ and $\o(\b)\in\pm\dBB$.
\nl
In {\it loc.cit.} it was conjectured that the signs in (a) are $+$.
The fact that the sign is $+$ for $\o\sha$  is proved in \cite{K94,4.3.2}. I thank H.Nakajima for pointing out to me that
\cite{K94,4.3.2} together with
\cite{L95, 4.14} imply that the sign for $\o$
is $+$. See also 3.16(a) for a more precise statement.

\subhead 3.9\endsubhead
For $\l\in X^+$ let $\dBB[\l]$ be the set of all
$\b\in\dBB\cap\dU[\ge\l]$ such that $\b$ acts on $\L_\l$ by a
nonzero map. By \cite{L93, 29.1.2, 29.1.3, 29.1.4} we have a partition
$\dBB=\sqc_{\l\in X^+}\dBB[\l]$. Note that for $\l\in X^+$, $\dU[\ge\l]$
(resp. $\dU[>\l]$) is the subspace of $\dU$ with basis
$\sqc_{\l'\in X^+;\l'\ge\l}\dBB[\l']$ (resp.
$\sqc_{\l'\in X^+;\l'>\l}\dBB[\l']$).

\subhead 3.10\endsubhead
Let $\l\in X^+$.

By \cite{L95, 4.4(a)}, for $b_1\in\BB(\l),b_2\in\BB(\l)$, there exists
a unique element $\b\in\dBB[\l]$ such that
$b_11_\l b_2^{*+}-\b\in\dU[>\l]$. We set $\b=\b_\l(b_1,b_2)$.

By \cite{L95, 4.4(b)},

(a) {\it the map $f:\sqc_{\l\in X^+}\BB(\l)\T\BB(\l)@>>>\dBB[\l]$ given by
$(\l,(b_1,b_2))\m\b_\l(b_1,b_2)$ is bijective.}

\proclaim{Lemma 3.11} Let $\l\in X^+$, $b_1\in\BB(\l),b_2\in\BB(\l)$. For
any $r\in\ZZ$ we set 
$$u_r:=1_\l b_2^{*+}b_1^-1_\l-v^r(b_1^-\et_\l,b_2^-\et_\l)_\l1_\l.$$
Then for some $r=r_{b_2,\l}\in\ZZ$ we have $u_r\in\dU[>\l].$
\endproclaim
Here we write $r=r_{b_2,\l}$ instead of: $r$ depending on $b_2,\l$
but not on $b_1$. 

The following proof is almost copied from \cite{L95, 4.7}.

Since $1_\l\in\dU[\ge\l]$ and $\dU[\ge\l]$ is a two-sided ideal of $\dU$,
we have $u_r\in\dU[\ge\l]$ and it is enough to show that for some
$r=r_{b_2,\l}\in\ZZ$, $u_r$ acts as $0$ on $\L_\l$. Since
$u_r=u_r1_\l$ and $1_\l\L_\l$ is the line spanned by $\et_\l$, it is
enough to show that $u_r\et_\l=0$ for some $r=r_{b_2,\l}\in\ZZ$.
Since $u_r=1_\l u_r$, we have
$u_r\et_\l=z_r\et_\l$ for some $z_r\in\QQ(v)$. We have
$$z_r(\et_\l,\et_\l)_\l=(u_r\et_\l,\et_\l)_\l.$$
By the definition of $(,)_\l$, for some $r_0=(r_0)_{b_2,\l}\in\ZZ$
we have
$$(1_\l b_2^{*+}b_1^-1_\l\et_\l,\et_\l)_\l=
(b_1^-\et_\l,v^{r_0}\sha(1_\l b_2^{*+})\et_\l)_\l=
(b_1^-\et_\l,v^{r_0}b_2^-\et_\l)_\l,$$
so that
$$z_{r_0}(\et_\l,\et_\l)_\l=(b_1^-\et_\l,b_2^-\et_\l)_\l)
(v^{r_0}-v^{r_0}(\et_\l,\et_\l)_\l).$$
Since $(\et_\l,\et_\l)_\l=1$ we see that $z_{r_0}=0$ so that
$u_{r_0}\et_\l=0$. The lemma is proved.

\subhead 3.12 \endsubhead
In the setup of Lemma 3.11 let $b_0\in\BB(\l)$. By Lemma 3.11
and its proof, we can find $r_0=(r_0)_{b_2,\l}\in\ZZ$ such that
$u_{r_0}\et_\l=0$; we then have $b_0^-u_{r_0}\et_\l=0$. We see that
$$b_0^-1_\l b_2^{*+}b_1^-1_\l\et_\l=
v^{r_0}(b_1^-\et_\l,b_2^-\et_\l)_\l b_0^-\et_\l.$$
We now replace $\l,b_0,b_1,b_2$ by $\l^!,\k(b_0),\k(b_1),\k(b_2)$.
(Recall that $\k(b_0),\k(b_1),\k(b_2)$ are in $\BB(\l^!)$). We see
that for any $\l\in X^+$ and any $b_0,b_1,b_2$ in $\BB(\l)$ we can
find $\ti r_0=(\ti r_0)_{b_2,\l}\in\ZZ$ such that
$$\k(b_0)^-1_{\l^!}\k(b_2)^{*+}\k(b_1)^-1_{\l^!}\et_{\l^!}=
v^{\ti r_0}(\k(b_1)^-\et_{\l^!},\k(b_2)^-\et_{\l^!})_{\l^!}
\k(b_0)^-\et_{\l^!}.$$
Using 3.6(b), we deduce
$$\k(b_0)^-1_{\l^!}\k(b_2)^{*+}\k(b_1)^-1_{\l^!}\et_{\l^!}= 
v^{\ti r_0}(b_1^+\x_{\l^!},b_2^+\x_{\l^!})_{\l^!}
b_0^+\x_{\l^!}.\tag a$$

\proclaim{Lemma 3.13} Let $\l\in X^+$,
$b_1\in\BB(\l),b_2\in\BB(\l)$.
Let $\d=(\x_{\l^!},\x_{\l^!})_{\l^!}$. For any $r\in\ZZ$ we set 
$$u'_r:=1_{-\l}b_2^{*-}b_1^+1_{-\l}
-v^r\d\i(b_1^+\x_{\l^!},b_2^+\x_{\l^!})_{\l^!}1_{-\l}.$$
Then for some $r=r_{b_2,\l}\in\ZZ$ we have $u'_r\in\dU[>\l^!].$
\endproclaim
The proof is similar to that of Lemma 3.11.

Since $1_{-\l}\in\dU[\ge\l^!]$ and $\dU[\ge\l^!]$ is a two-sided ideal
of $\dU$ we have $u'_r\in\dU[\ge\l^!]$ and it is enough to show that for
some $r=r_{b_2,\l}\in\ZZ$,
$u'_r$ acts as $0$ on $\L_{\l^!}$. Since $u'_r=u'_r1_{-\l}$
and $1_{-\l}\L_{\l^!}$ is the line spanned by $\x_{\l^!}$, it is enough
to show that $u'_r\x_{\l^!}=0$ for some $r=r_{b_2,\l}\in\ZZ$. Since
$u'_r=1_{-\l}u'_r$, we have $u'_r\x_{\l^!}=z'_r\x_{\l^!}$ for some
$z'_r\in\QQ(v)$. We have
$z'_r(\x_{\l^!},\x_{\l^!})_{\l^!}=(u'_r\x_{\l^!},\x_{\l^!})_{\l^!})$.
By the definition of $(,)_{\l^!}$, for some $r'_0=(r'_0)_{b_2,\l}\in\ZZ$
we have
$$(1_{-\l}b_2^{*-}b_1^+1_{-\l}\x_{\l^!},\x_{\l^!})_{\l^!}=
(b_1^+\x_{\l^!},v^{r'_0}\sha(1_{-\l}b_2^{*-})\x_{\l^!})_{\l^!}=
(b_1^+\x_{\l^!},v^{r'_0}b_2^+\x_{\l^!})_{\l^!},$$
so that
$$z'_{r'_0}(\x_{\l^!},\x_{\l^!})_{\l^!}=
(b_1^+\x_{\l^!},b_2^+\x_{\l^!})_{\l^!})
(v^{r'_0}-v^{r'_0}\d\i(\x_{\l^!},\x_{\l^!})_{\l^!}).$$
Since $(\x_{\l^!},\x_{\l^!})_{\l^!}=\d\ne0$ (see 3.7), we see that 
$z'_{r'_0}=0$ so that $u'_{r'_0}\x_{\l^!}=0$. The lemma is proved.

\subhead 3.14\endsubhead
In the setup of Lemma 3.13, let $b_0\in\BB(\l)$. By Lemma 3.13
and its proof we can find $r'_0=(r'_0)_{b_2,\l}\in\ZZ$ such that
$u'_{r'_0}\x_{\l^!}=0$; we then have
$(b_0^+1_{-\l})u'_{r'_0}\x_{\l^!}=0$. We see that
$$b_0^+1_{-\l}b_2^{*-}b_1^+1_{-\l}\x_{\l^!}=
v^{r'_0}\d\i(b_1^+\x_{\l^!},b_2^+\x_{\l^!})_{\l^!}b_0^+\x_{\l^!}.$$

Comparing with 3.12(a), we deduce
$$v^{-r'_0}\d b_0^+1_{-\l} b_2^{*-}b_1^+1_{-\l}\x_{\l^!}=
v^{-\ti r_0}\k(b_0)^-1_{\l^!}\k(b_2)^{*+}\k(b_1)^-1_{\l^!}\et_{\l^!}.\tag a$$

\subhead 3.15\endsubhead
Let $\l\in X^+$, $b_0\in\BB(\l)$, $b_2\in\BB(\l)$.
Since $1_{-\l}\in\dU[\ge\l^!]$, $1_{\l^!}\in\dU[\ge\l^!]$ 
and $\dU[\ge\l^!]$ is a two-sided ideal of $\dU$ we have
$$b_0^+1_{-\l}b_2^{*-}\in\dU[\ge\l^!],$$
$$\k(b_0)^-1_{\l^!}\k(b_2)^{*+}\in\dU[\ge\l^!].$$
Let $\mu=b_0^+1_{-\l}b_2^{*-}-C\k(b_0)^-1_{\l^!}\k(b_2)^{*+}$
where $C=v^{r'_0-\ti r_0}\d\i$ with $\ti r_0,r'_0$ as in 3.12,
3.13 and $\d$ is as in 3.13. We show:
$$\mu\in\dU[>\l^!].\tag a$$
It is enough to show that $\mu$ acts as zero on $\L_{\l^!}$ or that
$\mu s=0$ for any $s\in B_{\l^!}$ that is, for any $s$ of the form
$s=b_1^+\x_{\l^!}=\k(b_1)^-\et_{\l^!}$ with $b_1\in\BB(\l)$. Thus, it is
enough to show that
$$b_0^+1_{-\l}b_2^{*-}b_1^+\x_{\l^!}
-C\k(b_0)^-1_{\l^!}\k(b_2)^{*+}\k(b_1)^-\et_{\l^!}=0$$
for any $b_1\in\BB(\l)$. This clearly follows from 3.14(a).

We have $b_0^+1_{-\l}b_2^{*-}=\o(b_0^-1_\l b_2^{*+})$,
$b_0^-1_\l b_2^{*+}=\b_\l(b_0,b_2)+\g$
and
$$\k(b_0)^-1_{\l^!}\k(b_2)^{*+}=\b_{\l^!}(\k(b_0),\k(b_2))+\g'$$
where $\g\in\dU[>\l],\g'\in\dU[>\l^!]$ so that (a) implies
$$\o(\b_\l(b_0,b_2)+\g)-C(\b_{\l^!}(\k(b_0),\k(b_2))+\g')\in\dU_{>\l^!}.$$

From \cite{L93, 29.3.1} we see that $\o(\dU[\ge\l])\sub\dU[\ge\l^!]$ and
$\o(\dU[>\l])\sub\dU[>\l^!]$. Thus $\o(\g)\in\dU[>\l^!]$. We see that
$$\o(\b_\l(b_0,b_2))-C\b_{\l^!}(\k(b_0),\k(b_2))\in\dU[>\l^!].$$

By 3.8(a) we have $\o(\b_\l(b_0,b_2))=\e\b'$ where $\e\in\{1,-1\}$,
$\b'\in\dBB$ is necessarily in $\dBB[\ge\l^!]$ and we have
$$\e\b'-C\b_{\l^!}(\k(b_0),\k(b_2))\in\dU[>\l^!].\tag b$$
If $\b'\in\dBB[>\l^!]$, then we have
$$C\b_{\l^!}(\k(b_0),\k(b_2))\in\dU[>\l^!],$$
contradicting $$\b_{\l^!}(\k(b_0),\k(b_2))\in\dBB[\l^!],C\ne0$$
Thus, we have $\b'\in\dBB[\l^!]$ so that (b) implies
$$\b'=\e C\b_{\l^!}(\k(b_0),\k(b_2)).$$
It follows that $\e C=1$ that is, $\d=\e v^{r'_0-\ti r_0}$. Since
$\d\in1+v\i\QQ[[v\i]]$ (see 3.7), we see that
$$\d=1,\e=1,\ti r_0=r'_0,\tag c$$
so that $C=1$. Thus we have the following result.

\proclaim{Proposition 3.16} (a) Let
$\l\in X^+,b_0\in\BB(\l),b_2\in\BB(\l)$. We have
$$\o(\b_\l(b_0,b_2))=\b_{\l^!}(\k(b_0),\k(b_2)).$$
In particular, $\o(\dBB(\l))=\dBB(\l^!)$.

(b) We have $(\x_\l,\x_\l)_\l=1$.
\endproclaim

\subhead 3.17\endsubhead
Let $\AA=\QQ[[v\i]]\cap\QQ(v)$.
Let $\bof_\AA$ be the $\AA$-submodule of $\bof$ with basis $\BB$.

For $i\in I$ let $\ti f_i:\bof@>>>\bof$, $\ti e_i:\bof@>>>\bof$ be the
linear maps defined in \cite{K91}.
From \cite{K91,L90b} we see that there are well
defined maps $\ph_i:\BB@>>>\BB$, $\e_i:\BB@>>>\BB\cup\{0\}$ such
that for any $b\in\BB$ we have $\ti f_i(b)=\ph_i(b)\mod v\i\bof_\AA$,
$\ti e_i(b)=\e_i(b)\mod v\i\bof_\AA$.
Recall that for $b,b'$ in $\BB$ we have $\ph_i(b)=b'$ if and only if
$\e_i(b')=b$.

For $\l\in\L^+$ let $\L_{\l,\AA}$ be the $\AA$-submodule of $\L_\l$
with basis $B_\l$.
 For $i\in I$ let $\tF_i:\L_\l@>>>\L_\l$,
$\tE_i:\L_\l@>>>\L_\l$ be the linear maps denoted by $\ti f_i$,
$\ti e_i$ in \cite{K91}.
From \cite{K91,L90b} we see that there are well
defined maps $\cf_i:B_\l@>>>B_\l\cup\{0\}$,
$\ce_i:B_\l@>>>B_\l\cup\{0\}$ such that for any $d\in B_\l$ we
have $\ti F_i(d)=\cf_i(d)\mod v\i\L_{\l,A}$
$\ti E_i(d)=\ce_i(d)\mod v\i\L_{\l,A}$.
Recall that for $d,d'$ in $B_\l$ we have $\cf_i(d)=d'$ if and only
if $\ce_i(d')=d$. From the definitions we have:

(a) If $b,b'$ in $\BB(\l)$, $i\in I$ satisfy
$b^-\et_l=\cf_i(b'{}^-\et_l)$ then $b=\ph_i(b')$.

\proclaim{Lemma 3.18} If $b,b'$ in $\BB(\l^!)$, $i\in I$ satisfy
$b^+\x_l=\cf_i(b'{}^+\x_l)$, then $b=\e_i(b')$.
\endproclaim
Consider the vector space isomorphism $\t\i:\L_{\l^!}@>>>\L_\l$, see
3.6(a). It induces an $A$-module isomorphism
$\L_{\l^!,A}@>>>L_{\l,A}$ (since $\t\i(B_{\l^!})=B_\l$); moreover,
we have $\t\i\ti F_i=\ti E_i\t\i:\L_{\l^!}@>>>\L_\l$.

Let $d\in B_{\l^!}$; we have $\t\i(d)\in B_\l$. Applying $\t\i$ to
$$\ti F_i(d)=\cf_i(d)\mod v\i\L_{\l^!,A}$$ we obtain
$$\ti E_i(\t\i(d))=\t\i(\cf_i(d))\mod v\i\L_{\l,A}.$$
We have also
$$\ti E_i(\t\i(d))=\ce_i(\t\i(d))\mod v\i\L_{\l,A}.$$
Thus
$$\t\i(\cf_i(d))=\ce_i(\t\i(d))\mod v\i\L_{\l,A}.$$
Since $\t\i(\cf_i(d)),\ce_i(\t\i(d))$ are in $B_\l\cup\{0\}$, it
follows that

(a) $\t\i(\cf_i(d))=\ce_i(\t\i(d))$.
\nl
Now let $d'=b^+\x_\l$, $d=b'{}^+\x_\l$. By assumption we have
$d'=\cf_i(d)$. We have $\t\i(\cf_i(d))=\t\i(d')=b^-\et_\l$,
$\t\i(d)=b'{}^-\et_\l$. Using now (a), we deduce
$b^-\et_\l=\ce_i(b'{}^-\et_\l)$, so that $b=\e_i(b')$.
The lemma is proved.

\head 4. The involution $\ph_K$ \endhead
\subhead 4.1\endsubhead 
Let $\kk$ be an algebraically closed field of characteristic $0$.
Let $G$ be a connected reductive group over $\kk$ with a fixed
pinning corresponding to the root datum dual to that in 3.1. Thus,
$G$ has a given maximal torus $T$, given Borel subgroups $B^+,B^-$
with intersection $T$ (with unipotent radicals $U^+,U^-$) and given
imbeddings of algebraic groups $x_i:\kk@>>>U^+$,
$y_i:\kk@>>>U^-$ ($i\in I$) with the usual properties (see for
example \cite{L94, 1.1}). Note
that $X$ (resp. $Y$) is now the group of homomorphisms of algebraic
groups $\Hom(\kk^*,T)$ (resp. $\Hom(T,\kk^*)$).

Let $\Om:G@>>>G$ be the (involutive) automorphism of $G$ such that
$\Om(x_i(a))=y_i(a)$, $\Om(y_i(a))=x_i(a)$ for $i\in I,a\in\kk$,
$\Om(t)=t\i$ for $t\in T$.

Let $g\m \Th(g)$ be the (involutive) antiautomorphism of $G$ such
that $\Th(x_i(a))=x_i(a)$, $\Th(y_i(a))=y_i(a)$ for
$i\in I,a\in\kk$,
$\Th(t)=t\i$ for $t\in T$. We have $\Th\Om=\Om\Th$.

For $i\in I$ we set $\ds_i=y_i(-1)x_i(1)y_i(-1)$ (an element in the
normalizer of $T$).

Let $\cb$ be the variety of Borel subgroups of $G$. Now $\Om$ (resp.
$\Th$) induces an involution $\cb@>>>\cb$ denoted again by $\Om$
(resp. $\Th$); now $\Om$ interchanges $B^+,B^-$, while $\Th$
preserves $B^+,B^-$.

\subhead 4.2\endsubhead
We now fix a semifield $K$ contained in $\kk^*$ with $+,\T$ induced
from $\kk$. Let $U^+_{>0}$ (resp. $U^-_{>0}$) be the totally
positive part of $U^+$ (resp. $U^-$) defined in terms of
$K$ in \cite{L94, 2.12}. In \cite{L94, 2.7} a family of bijections
$g^+_\ii:K^\nu@>>>U^+_{>0}$ indexed by the various $\ii\in\ci$ is
considered; it is also shown, using Bruhat decomposition, that each
of these maps in injective. Using \cite{L94, 2.5}, we see that these
maps define a (surjective) map $g^+:\cu_K@>>>U^+_{>0}$. Note that
for any $\ii$ we have $g^+_\ii=g^+f_\ii$ where
$f_\ii:K^\nu@>>>\cu_K$ is as in 1.3. Since $g^+_\ii$ is injective,
it follows that

(a) $f_\ii$ is injective.
\nl
As shown in 1.3, $f_\ii$ is surjective hence bijective so that
$g^+=g^+_\ii f_\ii\i$. Since $g^+_\ii$ is injective it follows that
$g^+$ is injective. But it is also surjective so that it is
bijective. Thus,

(b) $g^+:\cu_K@>>>U^+_{>0}$ is a bijection.
\nl
Similarly, in \cite{L94, 2.9} a family of bijections
$g^-_\ii:K^\nu@>>>U^-_{>0}$ indexed by the various $\ii\in\ci$
is considered. Now these maps define a map $g^-:\cu_K@>>>U^-_{>0}$.
As above we see that

(c) $g^-:\cu_K@>>>U^-_{>0}$ is a bijection.

For $A\in\cu_K$ we write $A^+=g^+(A)\in U^+_{>0}$,
$A^-=g^-(A)\in U^-_{>0}$; we have $A^+=\Om(A^-)$.

Let $T_{>0}$ be the subgroup of $T$ generated by
$\{\l(a);\l\in X=\Hom(\kk^*,T),a\in K\}$.

\proclaim{Lemma 4.3} Let $\ii=(i_1,i_2,\do,i_\nu)\in\ci$,
$\aa=(a_1,a_2,\do,a_\nu)\in K^\nu$. We set $i=i_\nu$. Let
$k\in[1,\nu]$. For any $b_*=(b_{k+1},b_{k+2},\do,b_\nu)\in
K^{\nu-k}$ there exists a unique
$b'_*=(b'_k,b'_{k+1},\do,b'_\nu)\in K^{\nu-k+1}$ such that
$$\align&x_{i_1}(a_1)x_{i_2}(a_2)\do x_{i_{k-1}}(a_{k-1})
x_{i_k}(a_k)y_{i_\nu}(b_\nu)y_{i_{\nu-1}}(b_{\nu-1})\do\\&
y_{i_{k+1}}(b_{k+1})\ds_{i_{k+1}}\ds_{i_{k+2}}\do\ds_{i_\nu}B^-\\&
=x_{i_1}(a_1)x_{i_2}(a_2)\do x_{i_{k-1}}(a_{k-1})y_{i_\nu}(b'_\nu)
\\&y_{i_{\nu-1}}(b'_{\nu-1})\do y_{i_k}(b'_k)
\ds_{i_k}\ds_{i_{k+1}}\do\ds_{i_\nu}B^-.\tag a\endalign$$
Moreover, we have

(b) $b'_\nu=b_\nu(1+a_kb_\nu)\i$ if $i_k=i$,
$b'_\nu=b_\nu$ if $i_k\ne i$.
\endproclaim
The proof is essentially a repetition of arguments in the proof in
\cite{L97, 3.2} (except for (b)); with this occasion we correct some
typos in that proof. We only have to prove existence; the uniqueness
is immediate. Now $x_{i_k}(a_k)y_{i_\nu}(b_\nu)$ is equal to
$y_{i_\nu}(b'_\nu)x_{i_k}(c)t$ where $b'_\nu$ is
as in (b), $c\in K$ and $t\in T_{>0}$.
Hence the left hand side of (a) is  of the form
$$\align&x_{i_1}(a_1)x_{i_2}(a_2)\do x_{i_{k-1}}(a_{k-1})
y_{i_\nu}(b'_\nu)x_{i_k}(c)y_{i_{\nu-1}}(b''_{\nu-1})\do\\&
y_{i_{k+1}}(b''_{k+1})\ds_{i_{k+1}}\ds_{i_{k+2}}\do\ds_{i_\nu}B^-\endalign$$
for some $(b''_{k+1},\do,b''_\nu)\in K^{\nu-k}$.
Using the usual commutation relations between $x_{i_k}(?)$
and $y_j(?)$ we see that the left hand side of (a) is  of the form
$$\align&x_{i_1}(a_1)x_{i_2}(a_2)\do x_{i_{k-1}}(a_{k-1})
y_{i_\nu}(b'_\nu)y_{i_{\nu-1}}(b''_{\nu-1})\do\\&
y_{i_{k+1}}(b''_{k+1})x_{i_k}(c')\ds_{i_{k+1}}\ds_{i_{k+2}}\do\ds_{i_\nu}B^-\endalign$$
where $c'\in K$. We now use that
$x_{i_k}(c')=y_{i_k}(1/c')\ds_{i_k}y_{i_k}(1/c')t'$
where $t'\in T_{>0}$. we see that the left hand side of (a) is  of
the form
$$\align&x_{i_1}(a_1)x_{i_2}(a_2)\do x_{i_{k-1}}(a_{k-1})
y_{i_\nu}(b'_\nu)y_{i_{\nu-1}}(b''_{\nu-1})\do\\&
y_{i_{k+1}}(b''_{k+1})y_{i_k}(1/c')\ds_{i_k}y_{i_k}(1/c')
\ds_{i_{k+1}}\ds_{i_{k+2}}\do\ds_{i_\nu}B^-\endalign$$
It remains to observe that

$y_{i_k}(1/c')\ds_{i_{k+1}}\ds_{i_{k+2}}\do\ds_{i_\nu}B^-
=\ds_{i_{k+1}}\ds_{i_{k+2}}\do\ds_{i_\nu}B^-$

since 
$s_{i_k}s_{i_{k+1}}s_{i_{k+2}}\do s_{i_\nu}$
is a reduced expression in $W$. This proves the lemma.

\proclaim{Lemma 4.4} Let $\ii=(i_1,i_2,\do,i_\nu)\in\ci$,
$\aa=(a_1,a_2,\do,a_\nu)\in K^\nu$. We set $i=i_\nu$. Let
$k\in[1,\nu]$ and let $(b_{k+1},b_{k+2},\do,b_\nu)\in K^{\nu-k}$.
Assume that $b_\nu=(\sum_{l\in[k+1,\nu];i_l=i}a_l)\i$.
Then there exists a unique
$$(b'_k,b'_{k+1},\do,b'_\nu)\in K^{\nu-k+1}$$
such that
$$\align&x_{i_1}(a_1)x_{i_2}(a_2)\do x_{i_k}(a_k)
y_{i_\nu}(b_\nu)y_{i_{\nu-1}}(b_{\nu-1})\\&\do y_{i_{k+1}}(b_{k+1})
\ds_{i_{k+1}}\ds_{i_{k+2}}\do\ds_{i_\nu}B^-\\&
=x_{i_1}(a_1)x_{i_2}(a_2)\do x_{i_{k-1}}(a_{k-1})y_{i_\nu}(b'_\nu)
y_{i_{\nu-1}}(b'_{\nu-1})\\&\do y_{i_k}(b_k)\ds_{i_k}\ds_{i_{k+1}}
\ds_{i_{k+2}}\do\ds_{i_\nu}B^-.\endalign$$

Moreover, we have $b'_\nu=(\sum_{l\in[k,\nu];i_l=i}a_l)\i$.
\endproclaim
Except for the last sentence, this is a special case of Lemma 4.3.
It remains to prove the formula for $b'_\nu$. If $i_k\ne i$, then
from 4.3 we have
$$b'_\nu=(\sum_{l\in[k+1,\nu];i_l=i}a_l)\i
=(\sum_{l\in[k,\nu];i_l=i}a_l)\i.$$

If $i_k=i$, then from 4.3 we have
$$b'_\nu=\fra{(\sum_{l\in[k+1,\nu];i_l=i}a_l)\i}
{1+a_k(\sum_{l\in[k,\nu];i_l=i}a_l)\i}
=(\sum_{l\in[k,\nu];i_l=i}a_l)\i.$$

The lemma is proved.

\proclaim{Proposition 4.5} Let $\ii=(i_1,i_2,\do,i_\nu)\in\ci$,
$\aa=(a_1,a_2,\do,a_\nu)\in K^\nu$. We set $i=i_\nu$. 
There exists a unique $(b_1,b_2,\do,b_\nu)\in K^\nu$ such that
$$x_{i_1}(a_1)x_{i_2}(a_2)\do x_{i_\nu}(a_\nu)B^-
=y_{i_\nu}(b_\nu)y_{i_{\nu-1}}(b_{\nu-1})\do y_{i_1}(b_1)
\ds_{i_1}\ds_{i_2}\do\ds_{i_\nu}B^-.$$
Moreover, we have $b_\nu=(\sum_{l\in[1,\nu];i_l=i}a_l)\i$.
\endproclaim

Again this is contained in \cite{L97, 3.2} except for the last
sentence. It follows by applying 4.4 repeatedly with
$k=\nu,\nu-1,\do,1$.

\subhead 4.6\endsubhead
 From 4.5 we see that the image of the (injective) map

(a) $\cu_K@>c^+>>\cb$, $A\m A^+B^-(A^+)\i$
\nl
is contained in the image of the (injective) map

(b) $\cu_K@>c^->>\cb$, $A\m A^-B^+(A^-)\i$.
\nl
Applying $\Om$ we see that the image of $c^-$ is contained in the
image of $c^+$ hence

(c) the image of $c^+$ is equal to the image of $c^-$.
\nl
We denote this image by $\cb_{K,>0}$ (a subset of $\cb$). 
Applying $\Th$ to (c) we see that

(d) the image of $\cu_K@>>>\cb$, $A\m(A^+)\i B^-A^+$ is equal to
the image of $\cu_K@>>>\cb$, $A\m(A^-)\i B^+A^-$ and that these two
maps are injective.

We denote this image by $\cb_{K,<0}$ (a subset of $\cb$ equal to
$\Th(\cb_{K,>0})$). We see that $\Om:\cb@>>>\cb$ restricts to an
involution $\Om'_K:\cb_{K,>0}@>>>\cb_{K,>0}$ and to an involution
$\Om_K:\cb_{K,<0}@>>>\cb_{K,<0}$.

We also see that there is a unique bijection $\ph'_K:\cu_K@>>>\cu_K$
such that $$A^+B^-(A^+)\i=\ph'_K(A)^-B^+(\ph'_K(A)^-)\i$$ for any
$A\in\cu_K$, that is, $\Om'_K c^-=c^-\ph'_K$. Since $\Om'_K{}^2=1$,
it follows that $\ph'_K{}^2=1$. Moreover we see that there is a
unique bijection $\ph_K:\cu_K@>>>\cu_K$ such that
$$(A^+)\i B^-A^+=(\ph_K(A)^-)\i B^+\ph_K(A)^-$$ for any $A\in\cu_K$;
we have $\ph_K=*\ph'_K*$ ($*$ as in 1.1) hence $\ph_K^2=1$.

\subhead 4.7\endsubhead
From the proof of 4.5 we see that $\ph_K$ (and also $\ph'_K$),
when expressed as a map $K^\nu@>>>K^\nu$, is given by 
rational functions which are quotients of nonzero polynomials in
several variables with all coefficients in $\NN$ (after
substituting
elements of $K$ for the variables); moreover these polynomials
are independent of $K$. 
Replacing the variables in these polynomials by elements of an
arbitrary semifield $K'$ we see that $\ph_{K'}:\cu_{K'}@>>>\cu_{K'}$
is defined for any semifield $K'$ (not necessarily contained in a
field). Then various properties of this map can be deduced from the
analogous properties in the case where $K$ is as in 4.2.

\subhead 4.8\endsubhead
We can identify $\cb_K$ (see 1.8) with $\cb_{K,<0}$ by
$(A,A')\m(A^+)\i B^-A^+$. 

\head 5. Connecting bases with objects over the semifield $\ZZ$
\endhead
\subhead 5.1\endsubhead
Let $\uu:\cu_\NN@>\si>>\BB$ be the bijection defined in
\cite{L90a}. It satisfies $\uu(\hh_0)=1$.
More generally, if $A\in\cu_\NN$ and $b=\uu(A)$ then
$$wt(b)=\sum_{h\in I}||A||_hh'.\tag a$$
See \cite{L90a, 2.9}. From \cite{L90, 2.11} we see that under
$\uu$, the restriction to $\cu_\NN$ of the involution $*\io=\io*$
(see 1.1, 1.6, with $K=\ZZ$) corresponds to the involution $*$ of
$\BB$ in 3.2.

\subhead 5.2\endsubhead
In \cite{L90b} it is shown that for $i\in I$ we have
$$\ph_i\uu(A)=\uu(T_{i,1}A)\text{ for all }A\in\cu_\NN.$$

From (a) we can deduce:

(b) if $A\in\cu_\NN$ is such that $\e_i\uu(A)\in\BB$, then
$T_{i,-1}A\in\cu_\NN$ and $\e_i\uu(A)=\uu(T_{i,-1}A)$.
\nl
Indeed, we have $\e_i\uu(A)=\uu(A_1)$ for some $A_1\in\cu_\NN$.
Then $\uu(A)=\ph_i\uu(A_1)$ so that by (a) we have $A=T_{i,1}A_1$
hence $A_1=T_{i,-1}A$.

\subhead 5.3\endsubhead
Until the end of 5.6 we fix $\l\in X^+$. We define
$p^\l=(p^\l_i)\in\NN^I$ by $p^\l_i=(i,\l)$ for $i\in I$. We show:
$$\sum_{h\in I}||\hh_{p^\l}||_hh'=\l+\l^!.\tag a$$
Using 2.3(a) we see that this is equivalent to the equality
$$\sum_{h\in I}\sum_{i\in I}((i,\l)+(i,\l^!))b_{ih}h'=\l+\l^!.$$
that is 
$$\sum_{h\in I}\sum_{i\in I}(i,\z)b_{ih}h'=\z\tag b$$
where $\z=\l+\l^!=\l-w_0(\l)\in\sum_{h\in I}h'\sub X$.
It is enough to show that the two sides of (b) have the same $(,)$
with any $j\in I$ (viewed as an element of $Y$) that is,
$$\sum_{h\in I}\sum_{i\in I}(i,\z)b_{ih}a_{hj}=(j,\z)$$
for $j\in I$.
This follows immediately from the definition of $b_{ih}$.

\subhead 5.4\endsubhead
According to \cite{L90a},

(a) $\uu$ restricts to a bijection
$\uu^\l:\cu_{\NN,p^\l}@>\si>>\BB(\l)$.

We have $\uu^\l(\hh_0)^-\et_\l=\et_\l$.

We show:
$$\uu^\l(\hh_{p^\l})^-\et_\l=\x_\l.\tag b$$

Let $d=\uu^\l(\hh_{p^\l})^-\et_\l$. Since $d\in B_\l$ and
$\x_\l$ is the unique element $b\in B_\l$ such that
$1_{w_0(\l)}b=b$, it is enough to show that $1_{w_0(\l)}d=d$. Since
$1_\l\et_\l=\et_\l$, we see from the definition of $d$ that
$1_{\l'}d=d$ where $\l'=\l-wt(\uu^\l(\hh_{p^\l})$. Using 5.1(a),
5.3(a) this can be rewritten as $\l'=\l-(\l+\l^!)=-\l^!=w_0(\l)$.
This proves (b).

\proclaim{Theorem 5.5} (i) For any $A\in\cu_{\NN,p^\l}$ we have

$\io S_{p^\l}\ph_\ZZ(A)\in\cu_{\NN,p^{\l^!}}$.

(ii) Define a bijection $\ti\k:\cu_{\NN,p^\l}@>>>\cu_{\NN,p^{\l^!}}$ by
$\uu_{\l^!}(\ti\k(A))=\k(\uu_\l(A))$ for any $A\in\cu_{\NN,p^\l}$.
For
any $A\in\cu_{\NN,p^\l}$ we have

$\ti\k(A)=\io S_{p^\l}\ph_\ZZ(A)$.
\endproclaim
Something close to this is proved in \cite{L97, 4.9}; but that     
proof contains some misprints; for this reason we reprove it
without referring to \cite{L97, 4.9}.

Let $A\in\cu_{\NN,p^\l}$. Let $d=\uu^\l(A)^-\et_\l\in B_\l$. From
\cite{K91} it is known that there exists a sequence
$i_1,i_2,\do,i_k$ in $I$ such that
$$d=\cf_{i_1}\cf_{i_2}\do\cf_{i_k}\et_\l.$$
The smallest such $k$ is
a number $f(A)\in\NN$. We argue by induction on $f(A)$. If $f(A)=0$
we have $\uu^\l(A)^-\et_\l=\et_\l$ hence $\uu^\l(A)=1$ and
$A=\hh_0$. Using 2.2(c),(b),(a), we have
$$\io S_{p^\l}\ph_\ZZ(\hh_0)=\io S_{p^\l}(\hh_0)=\io \hh_{p^\l}=
\hh_{p^{\l^!}}$$
which belongs to $\cu_{\NN,p^{\l^!}}$; thus (i) holds in this case. By
the proof of (i), proving (ii) for $A=\hh_0$ it is the same as
proving that
$$\k(\uu^\l(\hh_0))=\uu^{\l^!}(\hh_{\p^{\l^!}})$$

or, using 5.4(b) for $\l^!$ instead of $\l$, that
$\k(1)^-\et_{\l^!}=\x_{\l^!}$.
From the definition of $\k$ this is the same as
$1^+\x_{\l^!}=\x_{\l^!}$ which is obvious.
Thus (ii) is proved in our case.

We can now assume that $A$ is such that $f(A)\ge1$ and that the
result is known when $A$ is replaced by any $A'\in\cu_{\NN,p^\l}$ such
that $f(A')<f(A)$. By the definition of $f(A)$ we can find
$A'\in\cu_{\NN,p^\l}$ such that $f(A')=f(A)-1$ and
$$\uu^\l(A)^-\et_\l=\cf_i(\uu^\l(A')^-\et_\l)$$
for some $i\in I$.
From 3.17(a) we then have $\uu(A)=\ph_i(\uu(A'))$ hence, by 5.2, we
have $\uu(A)=\uu(T_{i,1}A')$ so that $A=T_{i,1}A'$.

Let $b=\k(\uu_\l(A))=\k(\ph_i(\uu(A'))$, $b'=\k(\uu(A'))$. We have
$$b^+\x_\l=\uu_\l(A)^-\et_\l=(\ph_i\uu(A'))^-\et_\l=
\cf_i(\uu(A')^-\et_\l)=\cf_i(b'{}^+\x_\l).$$

Using 3.18 we see that $b=\e_i(b')$, that is

(a) $\k(\uu_\l(A))=\e_i\k(\uu_\l(A'))$.
\nl
By the induction hypothesis we have

(b) $\io S_{p^\l}\ph_\ZZ(A')\in\cu_{\NN,p^{\l^!}}$ and
$\ti\k(A')=\io S_{p^\l}\ph_\ZZ(A')$.
\nl
Using (a) we have 
$$\uu(\ti k(A))=\k(\uu_\l(A))=\e_i\k(\uu_\l(A'))
=\e_i\uu(\io S_{p^\l}\ph_\ZZ(A')).$$

(The third equality follows from (b).) Using 5.2(b) we see that
$$T_{i,-1}\io S_{p^\l}\ph_\ZZ(A')\in\cu_\NN$$
and

$$\e_i\uu(\io S_{p^\l}\ph_\ZZ(A'))
=\uu(T_{i,-1}\io S_{p^\l}\ph_\ZZ(A')).$$
Here the right hand side is equal to
$$\uu(\io S_{p^\l}\ph_\ZZ T_{i,1}(A'))=
\uu(\io S_{p^\l}\ph_\ZZ A).$$

(We have used 1.6(a),(b), 1.7(c).)

Thus we have $\uu(\ti k(A))=\uu(\io S_{p^\l}\ph_\ZZ A)$, so that
$\ti k(A)=\io S_{p^\l}\ph_\ZZ A$. We see that (i) and (ii) hold for
$A$. The theorem is proved.

\proclaim{Corollary 5.6} (i) For any $A\in\cu_{\NN,p^\l}$ we have

$S_{p^\l}\ph_\ZZ(A)\in\cu_{\NN,p^\l}$.

(ii) The map $\cu_{\NN,p^\l}@>>>\cu_{\NN,p^\l}$ given by
$A\m S_{p^\l}\ph_\ZZ(A)$ is a bijection.
\endproclaim
(i) follows from 5.5(i); (ii) follows from 5.5(ii).

\proclaim{Corollary 5.7} Let $p\in\NN^I$, $A\in\cu_{\NN,p}$. Then
$S_p\ph_\ZZ(A)\in\cu_{\NN,p}$. Moreover, the map
$\cu_{\NN,p}@>>>\cu_{\NN,p}$ given by $A\m S_p\ph_\ZZ(A)$ is a
bijection.
\endproclaim
We replace $Y,X,(,)$ by $Y',X'=\Hom(Y',\ZZ),(,)'$ where
$Y'=\sum_i\ZZ i\sub Y$ and $(,)'$ is the obvious pairing. Define
$i'\in X'$ by $(j,i')'=a_{ij}$. We apply 5.6 to this new root datum
and to $\l$ replaced by $\l'\in X'$ given by $(i,\l')=p_i$. Then 5.7
follows. 

\proclaim{Corollary 5.8} For any $p\in\NN^I$ we have
$$\cu_{\NN,p}=\{A\in\cu_\NN;S_p\ph_\ZZ(A)\in\cu_\NN\}.$$
Hence $(A,A')\m A$ is a bijection $\cb_{\NN,p}@>\si>>\cu_{\NN,p}$.
\endproclaim
The first sentence follows immediately from 5.7, as shown in
\cite{L21, no.12}. The
second sentence follows from the first sentence using 2.4(a).

\proclaim{Corollary 5.9} Let $\l\in X^+$. The map
$\cb_{\NN,p^\l}@>>>\BB(\l)$ given by $(A,A')\m\uu(A)$ is a
bijection.
\endproclaim
This follows from 5.4(a) and 5.8.

\subhead 5.10\endsubhead
Let
$$\align&\Xi=\sqc_{\l\in X^+}\cb_{\NN,p^\l}\T\cb_{\NN,p^\l}\\&
=\sqc_{\l\in X^+}\{(B_1,B_2)\in\cb_\ZZ\T\cb_\ZZ;B_1\in\cb^+_\NN,
B_2\in\cb^+_\NN,\un S_{p^{-\l}}B_1\in\cb^-_\NN,\\&
\un S_{p^{-\l}}B_2\in\cb^-_\NN\}.\endalign$$
We define $\dot\uu:\Xi@>>>\dBB$ by
$\dot\uu(B_1,B_2,\l)=b_\l(\uu(A_1),\uu(A_2))$
for any $\l\in X^+$ and any $B_1=(A_1,A'_1)\in\cb_{\NN,p^\l}$,
$B_2=(A_2,A'_2)\in\cb_{\NN,p^\l}$.

We define $\ti\o:\Xi@>>>\Xi$ by

$\ti\o(B_1,B_2,\l)=(\un\io\un\ph_\ZZ(B_1),\un\io\un\ph_\ZZ(B_2),\l^!)$.

We define $\ti\sha:\Xi@>>>\Xi$ by $\ti\sha(B_1,B_2,\l)=(B_2,B_1,\l)$.

\proclaim{Theorem 5.11}(a) $\dot\uu$ is a bijection.

(b) We have $\o\dot\uu=\dot\uu\ti\o$.

(c) We have $\sha\dot\uu=\dot\uu\ti\sha$.
\endproclaim
(a) follows from 3.10(a) and 5.9; (b) follows from 3.16, 5.5 and
5.8; (c) follows from \cite{L95, 4.14(a)}.

\subhead 5.12\endsubhead
The group $X$ acts on $\cb_\ZZ\T\cb_\ZZ$ by
$\l:(B,\tB)\m(\un S_{p^\l}B,\un S_{p^\l}\tB)$.

We now assume that the root datum in 3.1 is of simply connected
type, that is, the map $X@>>>\ZZ^I$, $\l\m p^\l$ is a bijection.
Then the $X$-action above is free. Let
$$\align&\Xi'=\{(B_1,B_2,\tB_1,\tB_2)\in
\cb^+_\NN\T\cb^+_\NN\T\cb^-_\NN\T\cb^-_\NN;\\&
(B_1,B_2)\text{ is in the $X$-orbit of }(\tB_1,\tB_2)\}.\endalign$$
We define $\Xi@>>>\Xi'$ by
$$(B_1,B_2,\l)\m(B_1,B_2,\un S_{p^{-\l}}B_1,\un S_{p^{-\l}}B_2).
\tag a$$
We show:

(b) The map (a) is a bijection.
\nl
The injectivity follows from the freeness of the $X$-action.
We show that (a) is surjective. It is enough to show that:

(c) if $B_1\in\cb^+_\NN$, $\tB_1\in\cb^-_\NN$, $\l\in X$ are such
that $\tB_1=\un S_{p^{-\l}}B_1$, then $\l\in X^+$.
\nl
This follows from 2.5(a).

Using (b), we see that in our case, $\Xi'$ can be viewed as an
indexing set for $\dBB$.

\head Appendix\endhead
\subhead A.1 \endsubhead
Let $p=(p_i)\in\NN^I$. Then the set
$\cu_{\NN,p}$ is finite since it is
in bijection with the finite set $B_\l$ for some $\l\in X^+$
(assuming that the root datum is of simply connected type).
But one would like to have a proof of the finiteness
of $\cu_{\NN,p}$ independent of the theory of canonical bases.
Such a proof will be given in this appendix.

Let $n\in\NN$. Let $\cu_\NN^n$ be the set of all $A\in\cu_\NN$
such that for any $(\ii,\aa)\in A$ we have $a_1\le n$.

\proclaim {Lemma A.2} Let $A\in\cu_\NN^n$. For any
$(\ii,\aa)\in A$ we have $a_k\le 2^{k-1}n$ for $k\in[1,\nu]$.
\endproclaim
We argue by induction on $k$. For $k=1$ the result is clear. Now
assume that $k\ge2$. Let $(\ii,\aa)\in A$.
We set $a_{k-1}=a,a_k=b$, $i_{k-1}=i,i_k=j$.
If $a_{ij}=0$ then
$$((i_1,\do,i_{k-2},j,i,\do),(a_1,\do,a_{k-2},b,a,\do))\in A$$
and by the induction hypothesis we have $b\le 2^{k-2}n$ hence
$b\le 2^{k-1}n$ as desired. Thus we can assume that $a_{ij}=-1$.

Case 1. Assume that $s_is_{i_{k+1}}\do s_{i_\nu}$ is not a
reduced expression; then
$$s_{i_{k+1}}\do s_{i_\nu}=s_is_{j_{k+2}}\do s_{j_\nu}$$ (reduced
expression) for some $j_{k+2},\do,j_\nu$ in $I$ so that for
some $b_{k+1},\do,b_\nu$ we have
$$((i_1,\do,i_{k-2},i,j,i,j_{k+2},\do,j_\nu),
(a_1,...a_{k-2},a,b,b_{k+1},...,b_\nu))\in A.$$
Here we can replace $i,j,i$ by $j,i,j$ and $a,b,b_{k+1}$ by
$$b+b_{k+1}-\min(a,b_{k+1}),\min(a,b_{k+1}),a+b-\min(a,b_{k+1}).$$
From the induction hypothesis we have
$b+b_{k+1}-\min(a,b_{k+1})\le 2^{k-2}n$.
If $a\le b_{k+1}$, then  $b+b_{k+1}-a\le 2^{k-2}n$
hence $$b\le b+b_{k+1}\le a+2^{k-2}n\le2^{k-2}n+2^{k-2}n=2^{k-1}n.$$
If $b_{k+1}\le a$ then $b\le 2^{k-2}n$ hence $b\le 2^{k-1}n$ as
desired.

Case 2. We can now assume that
$s_is_{i_{k+1}}\do s_{i_\nu}$ is a reduced expression. Then
setting $y=s_{i_{k+1}}\do s_{i_\nu}$ we have that $|s_iy|>|y|$,
$|s_jy|>|y|$, hence $|s_js_is_jy|=|y|+3$, hence can find
$u\in W$ such that $us_js_is_jy=w_0$, $|u|+|s_js_is_j|+|y|=\nu$.
Hence we can find $(\ii',\aa')\in A$ with

$(i'_{k-2},i'_{k-1},i'_k)=(j,i,j)$,

$(a'_{k-2},a'_{k-1},a'_k)=(a'_{k-2},a,b)$.
\nl
Here we can replace $j,i,j$ by $i,j,i$
and $a'_{k-2},a,b$ by
$$a+b-\min(a'_{k-2},b),\min(a'_{k-2},b),
a'_{k-2}+a-\min(a'_{k-2},b).$$

By the induction hypothesis we have
$a'_{k-2}\le 2^{k-3}n$, $a+b-\min(a'_{k-2},b)\le 2^{k-3}n$.
If $b\le a'_{k-2}$ then $b\le 2^{k-3}n$ hence $b\le 2^{k-1}n$.
If $a'_{k-2}\le b$  then $a+b-a'_{k-2}\le 2^{k-3}n$
hence $$b\le a+b\le a'_{k-2}+2^{k-3}n\le2^{k-3}n+2^{k-3}n=2^{k-2}n$$
hence $b\le 2^{k-1}n$.
This completes the induction step. The lemma is proved.

\proclaim{Proposition A.3} $\cu_\NN^n$ is a finite set.
\endproclaim
Let $\ii\in\ci$. Using 1.3(a) it is enough to show the finiteness of
the set of all $\aa\in\NN^\nu$ such that the equivalence class of
$(\ii,\aa)$ is in $\cu_\NN^n$. By the lemma, the number of elements
in this set is $\le(N+1)(2N+1)(4N+1)...(2^{\nu-1}N+1)$. The
proposition is proved.

\subhead A.4\endsubhead
We now choose $n$ such that $n\ge p_i$ for all $i\in I$.
Clarly, $\cu_{\NN,p}$ is contained in the image of $\cu_\NN^n$
under $A\m A^*$. Using A.3, we deduce that $\cu_{\NN,p}$ is finite.

\enddocument